\documentstyle[amscd,amssymb,verbatim,12pt]{amsart}
\newcommand{\ad}{\operatorname{ad}}
\newcommand{\del}{\partial}
\newcommand{\ssl}{\frak{sl}}
\newcommand{\ch}{\operatorname{ch}}
\newcommand{\CH}{\operatorname{CH}}
\renewcommand{\mod}{\operatorname{mod}}

\newcommand{\OO}{{\cal O}}
\newcommand{\Sym}{\operatorname{Sym}}

\newcommand{\SL}{\operatorname{SL}}

\renewcommand{\P}{{\Bbb P}}

\newcommand{\si}{\sigma}

\newcommand{\de}{\delta}

\renewcommand{\ker}{\operatorname{ker}}
\newcommand{\im}{\operatorname{im}}

\numberwithin{equation}{section}

\newtheorem{thm}{Theorem}[section]
\newtheorem{prop}[thm]{Proposition}
\newtheorem{lem}[thm]{Lemma}
\newtheorem{cor}[thm]{Corollary}
\newenvironment{rem}{\vspace{3mm}\noindent
{\bf Remark.}}{\vspace{3mm}}

\newenvironment{rems}{\vspace{3mm}
\noindent {\bf Remarks.}}{\vspace{3mm}}

\newcommand{\Pf}{\noindent {\it Proof}}

\newcommand{\ra}{\rightarrow}
\newcommand{\A}{{\Bbb A}}

\newcommand{\PP}{{\cal P}}

\newcommand{\SS}{{\cal S}}
\newcommand{\LL}{{\cal L}}

\newcommand{\Jac}{\operatorname{Jac}}
\renewcommand{\a}{\alpha}
\renewcommand{\b}{\beta}

\newcommand{\De}{\Delta}
\newcommand{\la}{\lambda}

\newcommand{\Z}{{\Bbb Z}}
\newcommand{\Q}{{\Bbb Q}}

\newcommand{\wt}{\widetilde}

\newcommand{\alg}{\operatorname{alg}}
\newcommand{\sub}{\subset}
\newcommand{\ed}{\qed\vspace{3mm}}

\title{Universal algebraic equivalences between tautological cycles on 
Jacobians of curves}
\author{A. Polishchuk}
\thanks{Supported in part by NSF grant}

\begin{document}
\begin{abstract} We present a collection of algebraic equivalences between
tautological cycles on the Jacobian $J$ of a curve, i.e., cycles in the
subring of the Chow ring of $J$ generated by the classes of
certain standard subvarieties of $J$. 
These equivalences are universal in the sense that they hold for all curves
of given genus. We show also that they are compatible with the action of
the Fourier transform on tautological cycles and compute this action explicitly.
\end{abstract}
\maketitle

\bigskip

\centerline{\sc Introduction}

\medskip

Let $J$ be the Jacobian of a smooth projective complex curve $C$ of
genus $g\ge 2$. For every $d$, $0\le d\le g$, consider the morphism
$\si_d:\Sym^d C\ra J:D\mapsto \OO_C(D-dp)$,
where $p$ is a fixed point on $C$ (for $d=0$ this is an embedding of the
neutral point into $J$). It is well-known that $\si_d$
is birational onto its image. Let us denote by
$w_d=[\si_{g-d}(\Sym^{g-d})]\in\CH^d(J)$, where $0\le d\le g$, the 
corresponding classes in the Chow ring $\CH(J)$ of $J$. Following Beauville who
studied the subring in $\CH(J)$ generated by these classes (see \cite{Bmain})
we call an element of this subring a {\it tautological cycle} on $J$.
Note that $w_0=1$, while $w_1$
is the class of the theta divisor on $J$.  
The Poincar\'e formula
states that $w_d$ is homologically equivalent to 
$\frac{w_1^d}{d!}$. However, it is known that in
general this formula fails to hold modulo algebraic equivalence
(although it does hold for hyperelliptic curves, see \cite{Col}).
More precisely, Ceresa in \cite{C} has shown that if $C$ is generic 
of genus $g\ge 3$ then
$w_d$ is not algebraically equivalent to $[-1]^*w_d$ for $1<d\le g-1$,
where $[-1]^*$ is the involution of the Chow ring induced
by the inversion morphism $[-1]:J\ra J$.
This raises the problem of finding universal
polynomial relations between the classes $w_d$ that hold modulo algebraic
equivalence. In this paper we derive a number of such relations. We do not
know whether our set of relations is complete for a generic curve.
However, we show that these relations are in some sense consistent 
with the action of the Fourier transform on $\CH(J)$.

Let us set $p_k=-N^k(w)\in\CH^k(J)$ for $k\ge 1$, 
where $N^k(w)$ are the Newton polynomials on the classes $w_1,\ldots,w_g$:
$$N^k(w)=\frac{1}{k!}\sum_{i=1}^g\la_i^k,$$
where $\la_1,\ldots,\la_g$ are roots of the equation
$\la^g-w_1\la^{g-1}+\ldots+(-1)^gw_g=0$.
For example, $p_1=-w_1$, $p_2=w_2-w_1^2/2$, $p_3=w_2w_1/2+w_3/2-w_1^3/6$.
From the Poincar\'e formula it is easy to see that the classes
$p_n$ for $n>1$ are homologically trivial. Let us denote by
$\CH(J)_{\Q}/(\alg)$ the quotient of the Chow ring of $J$ with
rational coefficients modulo the ideal of classes algebraically equivalent
to $0$. 
It is also convenient to have a notation for
the divided powers of $p_1$: $p_1^{[d]}:=p_1^d/d!$.
Our main result is the following collection of relations in 
$\CH(J)_{\Q}/(\alg)$.

\begin{thm}\label{mainthm} 
(i) Let us define the differential operator $D$ acting on polynomials
in infinitely many variables $x_1,x_2,\ldots$:
$$D=-g\del_1+\frac{1}{2}
\sum_{m,n\ge 1}{m+n\choose m}x_{m+n-1}\del_m\del_n,$$
where $\del_i=d/dx_i$. 
Then for every polynomial $F$ of the form
$$F(x_1,x_2,\ldots)=D^d(x_1^{m_1}\ldots x_k^{m_k}),$$
where $m_1+2m_2+\ldots+km_k=g$, $m_1<g$ and $d\ge 0$, one has
$$F(p_1,p_2,\ldots)=0.$$
in $\CH^{g-d}(J)_{\Q}/(\alg)$.

\noindent
(ii) Here is another description of the same collection of relations.
For every $k\ge 1$, every $n_1,\ldots,n_k$ such that $n_i>1$, and every $d$
such that $0\le d\le k-1$, one has
\begin{equation}\label{mainid}
\sum_{[1,k]=I_1\sqcup I_2\sqcup\ldots\sqcup I_m}
{m-1\choose d+m-k}b(I_1)\ldots b(I_m)
p_1^{[g-d-m+k-\sum_{i=1}^k n_i]}p_{d(I_1)}\ldots p_{d(I_m)}=0
\end{equation}
in $\CH^{g-d}(J)_{\Q}/(\alg)$,
where the summation is over all partitions of the set $[1,k]=\{1,\ldots,k\}$
into the disjoint union of nonempty subsets $I_1,\ldots,I_m$
such that $-d+k\le m\le g-d+k-\sum_{i=1}^k n_i$ (two partitions
differing only by the ordering of the parts are considered to be the same);
for a subset $I=\{i_1,\ldots,i_s\}\sub[1,k]$ we denote 
$$b(I)=\frac{(n_{i_1}+\ldots+n_{i_s})!}{n_{i_1}!\ldots n_{i_s}!},$$
$$d(I)=n_{i_1}+\ldots+n_{i_s}-s+1.$$
\end{thm}  

Let us point out some corollaries of these relations.

\begin{cor}\label{CVGcor} The class $p_n$ is algebraically equivalent to
$0$ for $n\ge g/2+1$.
\end{cor}

In \cite{CG} Colombo and Van Geemen proved 
that $p_n$ is algebraically equivalent to $0$ for $n\ge d$, 
where $d$ is the minimal degree of a nonconstant
morphism from $C$ to $\P^1$ (see \cite{Bmain}, Proposition 4.1). 
For generic curve this is equivalent to the statement of the above
corollary.

\begin{cor}\label{lincombcor}
Every tautological cycle in $\CH^{g-d}(J)_{\Q}$ is
algebraically equivalent to a linear combination
of the classes $p_1^{g-d-\sum_{i=1}^k n_i}p_{n_1}\ldots p_{n_k}$, where
$0\le k\le d$ and $n_i>1$ for all $i=1,\ldots,k$. More precisely, 
for arbitrary $n_1,\ldots,n_k$ such that $n_i>1$ and for $d$ such that
$0\le d\le k-1$ one has
\begin{equation}\label{niceeq}
\begin{array}{l}
p_1^{[g-d-\sum_{i=1}^k n_i]}p_{n_1}\ldots p_{n_k}=\\
\sum_{j=1}^d(-1)^{k+d}{k-1-j\choose d-j}
\sum_{[1,k]=I_1\sqcup I_2\sqcup\ldots\sqcup I_j}b(I_1)\ldots b(I_j)
p_1^{[g-d-j+k-\sum_{i=1}^k n_i]}p_{d(I_1)}\ldots p_{d(I_j)}.
\end{array}
\end{equation}
in $\CH^{g-1}(J)_{\Q}/(\alg)$.
\end{cor}

For example, in the case $d=1$ the above corollary states that
for $n_1,\ldots,n_k$ such that $n_i>1$ and
$\sum_{i=1}^k n_i\le g-1$ one has
\begin{equation}\label{niceeq1}
p_1^{[g-1-\sum_{i=1}^k n_i]}p_{n_1}\ldots p_{n_k}=(-1)^{k-1}\frac{(\sum_{i=1}^k n_i)!}
{n_1!\ldots n_k!}p_1^{[g+k-2-\sum_{i=1}^k n_i]}p_{1-k+\sum_{i=1}^k n_i}
\end{equation}
in $\CH^{g-1}(J)_{\Q}/(\alg)$.

The main tool in the proof of Theorem \ref{mainthm}
is the Fourier transform $S:\CH(J)\ra \CH(J)$
introduced and studied by Beauville in \cite{B1}, \cite{B2} and \cite{Bmain}.
Theorem \ref{mainthm} is closely related to the
following explicit formula for the induced action of $S$ on
tautological cycles in $\CH(J)_{\Q}/(\alg)$.

\begin{thm}\label{fourthm}
For every $n_1,\ldots,n_k$ such that $n_i>1$ and every $n\ge 0$ one has
\begin{equation}\label{foureq}
(-1)^nS(p_1^{[n]}p_{n_1}\ldots p_{n_k})=
\sum_{[1,k]=I_1\sqcup I_2\sqcup\ldots\sqcup I_m}
b(I_1)\ldots b(I_m)
p_1^{[g-n-m-\sum_{i=1}^k n_i]}p_{d(I_1)}\ldots p_{d(I_m)}
\end{equation}
in $\CH(J)_{\Q}/(\alg)$,
where the summation is similar to the one in equation (\ref{mainid}).
\end{thm}


Note that the expression in the identity (\ref{mainid}) is 
equal to the RHS of equation (\ref{foureq}) in the case
$n=-1$, $d=k-1$. More interesting observation is the 
formal consistency between Theorems \ref{mainthm} and \ref{fourthm}
proved in part (iii) of the next theorem.

\begin{thm}\label{formalthm} For a fixed $g\ge 2$ let us denote
by $R^{\Jac}_g$ the quotient-space of
$\Q[x_1,x_2,\ldots]$ by the linear span of all the polynomials
appearing in Theorem \ref{mainthm} together with all polynomials
of degree $>g$ where $\deg(x_i)=i$. In other words, 
$$R^{\Jac}_g=\Q[x_1,x_2,\ldots]/I_g$$
where the subspace 
$I_g$ is spanned by all polynomials $F$ such that $\deg F>g$ and
by polynomials of the form $D^d(x_1^{m_1}\ldots x_k^{m_k})$,
where $d\ge 0$, $m_1+2m_2+\ldots+km_k=g$, $m_1<g$. Let us denote
by $p_i$ the image of $x_i$ in $R^{\Jac}_g$.
Then 

\noindent
(i) $I_g=\cap_{n\ge 1}\im(D^n)$, where $\im(D^n)$ denotes the image
of the operator $D^n$ acting on $\Q[x_1,x_2,\ldots]$;

\noindent
(ii) $I_g$ is an ideal in $\Q[x_1,\ldots,x_g]$, so $R^{\Jac}_g$ has
a commutative ring structure;

\noindent
(iii) the formula
(\ref{foureq}) gives a well-defined operator $S$ on $R^{\Jac}_g$
such that
$$S^2(p_1^np_{n_1}\ldots p_{n_k})=(-1)^{k+\sum_{i=1}^k n_i}
p_1^np_{n_1}\ldots p_{n_k},$$
where $n_i>1$; 

\noindent
(iv) the operators
\begin{equation}
e(F)=x_1\cdot F,\ f=-D,\ h=-g+\sum_{n\ge 1}(n+1)x_n\del_n
\end{equation}
on $\Q[x_1,x_2,\ldots]$ 
define a representation of the Lie algebra $\ssl_2$.
Furthermore, they preserve the ideal $I_g$ and therefore define
a representation of $\ssl_2$ on $R^{\Jac}_g$.
The operators $S$, $e$, $f$ and $h$ on $R^{\Jac}_g$ 
satisfy the standard compatibilities:
$$SeS^{-1}=-f,\ SfS^{-1}=-e,\ ShS^{-1}=-h.$$
\end{thm}

Thus, the ring $R^{\Jac}_g$ can be equipped with the same special
structures as $\CH(J)_{\Q}/(\alg)$ (namely,
the Fourier transform, the $\ssl_2$-action and the Pontryagin product) 
and the natural homomorphism $R^{\Jac}_g\ra\CH(J)_{\Q}/(\alg)$ 
respects these structures.

The next result points to possible connections between 
the structures on $R^{\Jac}_g$ and those on the 
cohomology of Hilbert schemes $(\A^2)^{[n]}$ of points on the plane. 
For example, the differential operators $D_k$ appearing below
are very similar to the operators describing the action of the Chern
character of the tautological bundle on $H^*((\A^2)^{[n]},\Q)$
(see \cite{Lehn},\cite{LS}). 

\begin{thm}\label{diffthm} 
For every $k\ge 2$ let us consider the following differential
operator on $\Q[x_1,x_2,\ldots]$:
$$D_k=\frac{1}{k!}
\sum_{n_1,\ldots,n_k\ge 1}\frac{(n_1+\ldots+n_k)!}{n_1!\ldots n_k!}
x_{n_1+\ldots+n_k-1}\del_{n_1}\ldots\del_{n_k}.$$
Note that the operator $D$ considered above is $D=D_2-g\del_1$.
Then 
for every $k\ge 3$ we have $[D_k,D]=0$, so $D_k$ descends to an operator
on $R^{\Jac}_g$. Furthermore,
for every $k\ge 2$ one has 
\begin{equation}\label{diffpk}
S(p_k\cdot S^{-1}F)=D_{k+1}(F),
\end{equation}
where $F\in R^{\Jac}_g$.
\end{thm}

\begin{rems} 1. Applying the homomorphism 
$R^{\Jac}_g\to\CH(J)_{\Q}/(\alg)$ we derive that the relation (\ref{diffpk}) also
holds in $\CH(J)_{\Q}/(\alg)$ whenever $F$ is a tautological cycle.
In fact, this relation can also be proved geometrically using 
similar arguments as for the case $k=1$, which corresponds to equation
(\ref{Ufor}). 

\noindent
2. The Lie algebra generated by operators $D_k$ for 
$k\ge 2$ and by the operators of multiplication by $x_k$ can be easily 
described: it is isomorphic to the Lie algebra of (polynomial)
hamiltonian vector fields
on the plane vanishing at the origin (where the plane is equipped with the
standard symplectic form). Elsewhere we will show that the action of
this Lie algebra on $R^{\Jac}_g$ extends to an action on the entire
Chow ring (not just tautological cycles).
\end{rems}

We should point out that at present it is not known whether there
exists a curve for which one of the classes $p_n$ with $n\ge 3$
is actually nonzero. Ceresa's theorem in \cite{C} can be interpreted
as nonvanishing of $p_2$ for generic curve of genus $g\ge 3$. Using
Theorem \ref{fourthm} we find that $S(p_2)=p_1^{[g-3]}p_2$. This implies
that for generic curve of genus $g\ge 3$ one has
$p_1^{g-3}p_2\neq 0$. We conjecture that for a generic curve the map
$$R^{\Jac}_g\ra\CH(J)_{\Q}/(\alg)$$
is injective, i.e., the set of relations of Theorem \ref{mainthm} is complete.

\bigskip

\noindent
{\it Acknowledgment}. I am grateful to Arkady Vaintrob for valuable discussions
and to Alexander Postnikov for providing
the references on Hurwitz identity (equation (\ref{combid2})).

\noindent
{\it Notation.} We use the notation ${x\choose n}=x(x-1)\ldots (x-n+1)/n!$,
where $n\ge 0$ and $x$ is either a number or an operator. Thus, we have 
${x\choose 0}=1$ while for $n>0$ and 
for an integer $m$ we have ${m\choose n}=0$ iff $m\in\{0,1,\ldots,n-1\}$.
If $n<0$ then we set ${x\choose n}=0$. 

\section{Fourier transform of cycles: brief reminder}

Let $X$ be an abelian variety of dimension $g$, $\hat{X}$ the dual abelian
variety, $\PP$ the Poincar\'e line bundle on $X\times\hat{X}$.
Recall that the Fourier-Mukai transform is an equivalence
$\SS:D^b(X)\ra D^b(\hat{X})$ sending a complex of coherent sheaves $F$
to $Rp_{2*}(p_1^*F\otimes\PP)$, where $p_1$ and $p_2$ are projections
of the product $X\times\hat{X}$ to its factors (see \cite{M}). 
It induces the transform
on Grothendieck groups and therefore on Chow groups with rational coefficients.
Explicitly, the Fourier transform of Chow groups 
$$S:\CH(X)_{\Q}\ra\CH(\hat{X})_{\Q}$$
is given by the formula
$S(\a)=p_{2*}(p_1^*\a\cdot\ch(\PP))$, where $\ch$ denotes the Chern character. 
This operator was introduced and studied by Beauville in \cite{B1},
\cite{B2}. It is clear that it preserves the ideal of cycles algebraically
(resp. homologically) equivalent to $0$. In particular, it induces a 
well-defined map 
$$S:\CH(X)_{\Q}/(\alg)\ra\CH(\hat{X})_{\Q}/(\alg).$$
For every $n\in\Z$ let us denote by $[n]_*,[n]^*:\CH(X)\ra \CH(X)$ the
operators
on the Chow group given by the push-forward and the pull-back with
respect to the endomorphism $[n]:X\ra X$.
The main properties of the Fourier transform that we will use are:
\begin{eqnarray}\label{conv}
S^2=(-1)^g[-1]^*,\nonumber\\
S(\a*\b)=S(\a)\cdot S(\b),\nonumber\\
S\circ [n]_*=[n]^*\circ S,
\end{eqnarray}
where $\a*\b$ denotes the Pontryagin product 
of $\a$ and $\b$ (the push-forward of
$\a\times\b$ be the group law morphism $X\times X\ra X$).

Beauville proved in \cite{B2} that there is a direct sum
decomposition of the Chow groups
$$\CH^d(X)_{\Q}=\oplus\CH^d_s(X),$$
where $\CH^d_s(X)=\{x\in\CH^d(X)_{\Q}:
[n]^*x=n^{2d-s}x\ \text{for all}\ n\in\Z\}$. 
This decomposition is compatible with various operations on the Chow groups
in the following way:
\begin{eqnarray}\label{bigreq}
\CH^p_s\cdot\CH^q_t\sub \CH^{p+q}_{s+t},\nonumber\\
\CH^p_s *\CH^q_t\sub \CH^{p+q-g}_{s+t},\nonumber\\
S(\CH^p_s)=\CH^{g-p+s}_s.
\end{eqnarray}

\begin{rem} If $X$ is principally polarized then identifying
$X$ with $\hat{X}$ we can view the Fourier transform $S$ as an automorphism
of $\CH(X)_{\Q}$. Let $c_1(L)\in\CH^1(X)$ be the class of the
principal polarization. Then
there is a natural action of the Lie algebra $\ssl_2$ on
$\CH(X)_{\Q}$ such that $e(x)=c_1(L)x$, $f(x)=-S(c_1(L)S^{-1}x)$,
$h(x)=(2p-g-s)x$ for $x\in\CH^p_s$. Indeed, this can be deduced 
from the explicit formulae for the
algebraic action of the group $\SL_2$ on $CH(X)_{\Q}$ considered in
\cite{thesis}. The corresponding action of $\ssl_2$
on tautological cycles in the Jacobian will play a crucial role below.
\end{rem}

\section{Computations in the Chow ring}\label{compseq}

\subsection{}\label{compseq1}

Let $C$ be a (connected) smooth projective curve of genus $g\ge 2$, 
$J$ be its Jacobian.
We can identify the $J$ with its dual $\hat{J}$ 
and consider the Fourier transform $\SS$ (resp., $S$)
as an autoequivalence of $D^b(J)$ (resp., as an automorphism of $\CH(J)$).
More precisely, to define this transform
we use the line bundle $\LL=\OO_{J\times J}(p_1^{-1}(\Theta)+p_2^{-1}(\Theta)-m^{-1}(\Theta))$ 
as a kernel on $J\times J$, where $\Theta\sub J$ is a theta divisor, $p_i$ are projections
of $J\times J$ to $J$, $m:J\times J\ra J$ is the group law.
Hence, $S$ induces a $\Q$-linear operator on $\CH(J)_{\Q}/(\alg)$.
The assertion of the following important lemma is proved in Proposition 2.3
and Corollary 2.4 of \cite{Bmain}.

\begin{lem}\label{chernlem}
$$S(w_{g-1})=\sum_{n=1}^{g-1}p_n$$
in $\CH(J)_{\Q}/(\alg)$.
Moreover, for every $n$ one has $p_n\in \CH^n_{n-1}/(\alg)$.
\end{lem}

For every collection of integers $n_1,\ldots,n_k$ we set
$$w(n_1,\ldots,n_k):=[n_1]_*w_{g-1}*\ldots*[n_k]_*w_{g-1}.$$
The following lemma is essentially contained in section (3.3) of \cite{Bmain}.

\begin{lem}\label{geomlem}
One has 
\begin{align*}
&w_1\cdot w(n_1,\ldots,n_k)=\sum_{i=1}^k[gn_i^2+n_i(\sum_{j\neq i}n_j)]
w(n_1,\ldots,\widehat{n}_i,\ldots,n_k)-\\
&\sum_{i<j} n_in_j
w(n_i+n_j,n_1,\ldots,\widehat{n}_i,\ldots,\widehat{n}_j,\ldots,n_k)
\end{align*}
where a hat over a symbol means that it should be omitted.
\end{lem}

\Pf . By definition $w(n_1,\ldots,n_k)$ is the push-forward of
the fundamental class under the composite map
$$u:C^k\stackrel{(\si)^k}{\ra} J^k\stackrel{{\bf n}}{\ra} J^k\stackrel{m}{\ra} J,$$
where $\si:C\ra J$ is the standard embedding,
${\bf n}=[n_1]\times\ldots\times[n_k]:J^k\ra J^k$,
$m$ is the addition morphism.
Therefore,
$$w_1\cdot w(n_1,\ldots,n_k)=u_*u^*w_1.$$
Using theorem of the cube one can easily show (see section (3.3)
of \cite{Bmain}) that
$$u^*w_1=\sum_i n_i^2 q_i^*\si^*w_1-
\sum_{i<j}n_in_j(\de_{ij}-q_i^*[p]-q_j^*[p])$$
modulo algebraic equivalence,
where $p\in C$ is a point,
$q_i$ are projections of $C^k$ to $C$,
$\de_{ij}\sub C^k$ is the class of the partial diagonal
divisor given by the equation $x_i=x_j$.
Since $\si^*w_1$ has degree $g$, it is algebraically equivalent
to $g[p]$. Hence, we have
$$u^*w_1=\sum_i (gn_i^2+n_i\sum_{j\neq i}n_j)q_i^*[p]-
\sum_{i<j}n_in_j\de_{ij}$$
modulo algebraic equivalence,
which immediately implies the assertion.
\ed

\subsection{}\label{compsec2}

The idea of the proof of Theorem \ref{fourthm} is to use the following
identity proved in \cite{Bmain}, (1.7):
\begin{equation}\label{fourexpeq}
S[-1]^*x=\exp(-p_1)\cdot\left(\exp(p_1)*[\exp(-p_1)\cdot x]\right).
\end{equation}
Thus, if we want to find explicitly the action of the Fourier
transform on polynomials in $p_i$, it suffices to find the formula
for the action of the operator $x\mapsto\exp(p_1)*x$ on such polynomials.
The only geometric ingredients needed for this are
Lemmata \ref{chernlem} and \ref{geomlem}, the rest 
of the computation consists of formal manipulations.
Along the way we will derive the relations of Theorem \ref{mainthm}.

\begin{lem}\label{Swlem} 
One has
$$S(w(n_1,\ldots,n_k))=n_1\ldots n_kP(n_1)\cdot\ldots\cdot P(n_k),$$
where $P(t)=\sum_{i=1}^{g-1}p_it^i$.
\end{lem}

\Pf . This follows immediately from Lemma \ref{chernlem} combined
with (\ref{conv}).
\ed

Let us define the operator $U$ on the space $\CH(J)_{\Q}/(\alg)$
by the formula
$$U(x)=S(p_1\cdot S^{-1}(x)).$$
Since $p_1$ is invariant with respect to $[-1]^*$ we also have
$$U(x)=S^{-1}(p_1\cdot S(x)).$$
Let us denote $U^{[n]}=U^{n}/n!$ for $n\ge 0$.
The plan of the computations below is the following.
First, we will compute the action of $U$ on a monomial in
$p_i$'s. Next, we will compute the action of the operators
$U^{[n]}$ on such a monomial for all $n\ge 0$.
As a result, we will get an explicit formula for the operator
$x\mapsto\exp(p_1)*x$ on the tautological ring. Finally,
using equation (\ref{fourexpeq}) we will find the formula
for the Fourier transform.

\begin{lem}\label{Uconvlem}
For every $n\ge 0$ one has
$$U^{[n]}(x)=(-1)^{g-n}p_1^{[g-n]}*x.$$
\end{lem}

\Pf . It is well known that $S(\exp(-p_1))=\exp(p_1)$.
Therefore, $(-1)^{g-n}S(p_1^{[g-n]})=p_1^{[n]}$.
Hence,
$$(-1)^{g-n}S(p_1^{[g-n]}*x)=p_1^{[n]}\cdot S(x)=SU^{[n]}(x).$$
\ed

\begin{lem}\label{convlem1} 
One has the following identity of polynomials in $t_1,\ldots,t_k$:
\begin{align*}
&U(P(t_1)\ldots P(t_k))=
-\sum_{i=1}^k[gt_i+\sum_{j\neq i}t_j]
P(t_1)\ldots\widehat{P(t_i)}\ldots P(t_k)+\\
&\sum_{i<j} (t_i+t_j)P(t_i+t_j)P(t_1)\ldots
\widehat{P(t_i)}\ldots,\widehat{P(t_j)}\ldots P(t_k),
\end{align*}
where hats denote omitted symbols.
\end{lem}

\Pf . Since both sides of the identity are polynomials in $t_1,\ldots,t_k$
it suffices to prove it for $t_i=n_i\in\Z$.
Using Lemma \ref{Swlem} we get
$$U(P(n_1)\ldots P(n_k))=
S(p_1\cdot \frac{w(n_1,\ldots,n_k)}{n_1\ldots n_k}).$$
It remains to apply the formula of Lemma \ref{geomlem} to compute
the expression under Fourier transform and then use Lemma \ref{Swlem}
again.
\ed

Thus, we arrive to the following formula for the action of the operator $U$
on monomials in classes $p_1,\ldots,p_{g-1}$.

\begin{prop}\label{convprop} 
For $n\ge 0$ and $n_1\ldots,n_k$ such that $n_i>1$ one has
\begin{align*}
&U(p_1^{[n]}p_{n_1}\ldots p_{n_k})=
(-g+n-1+k+\sum_{i=1}^k n_i) p_1^{[n-1]}p_{n_1}\ldots p_{n_k}+\\
&\sum_{i<j}{n_i+n_j\choose n_i} 
p_1^{[n]}p_{n_i+n_j-1}p_{n_1}\ldots\widehat{p}_{n_i}\ldots
\widehat{p}_{n_j}\ldots p_{n_k}.
\end{align*}
\end{prop}

\Pf . The expression in the LHS is equal to $1/n!$ times the coefficient with
$t_1^{n_1}\ldots t_k^{n_k}t^n$ in
$U(P(t_1)\ldots P(t_k)P(t)^n)$.
Now the assertion is an easy consequence of Lemma \ref{convlem1}.
\ed

We can immediately recognize in the formula of the above proposition
the action of the differential operator $D$ appearing in 
Theorem \ref{mainthm}, so that for every polynomial $F(x_1,x_2,\ldots)$
we have 
\begin{equation}\label{Ufor}
U(F(p_1,p_2,\ldots))=DF(p_1,p_2,\ldots).
\end{equation}

It is convenient to separate in $D$ the part that does not contain $\del_1$.
Namely, let us set
$$\De=\frac{1}{2}\sum_{m,n\ge 2}
{m+n\choose m}x_{m+n-1}\del_m\del_n.$$
Then we have 
\begin{equation}\label{Dfor}
D=\del_1 H+\De
\end{equation}
where 
$$H=-g-1+x_1\del_1+\sum_{n\ge 2}(n+1)x_n\del_n,$$ 
so that
$$H(x_1^nx_{n_1}\ldots x_{n_k})=
(-g+n-1+k+\sum_{i=1}^k n_i) x_1^nx_{n_1}\ldots x_{n_k},$$
where $n_i>1$ for $i=1,\ldots,k$. 

Let us set $\De^{[n]}=\De^n/n!$ for $n\ge 0$.

\begin{lem}\label{Upowerlem} For every polynomial $F(x_1,x_2,\ldots)$ and every
$n\ge 0$ one has
$$U^{[n]}(F(p_1,p_2,\ldots))=
\sum_{i=0}^n \del_1^{n-i}\De^{[i]}{H-i\choose n-i}F(p_1,p_2,\ldots).
$$
in $\CH(J)_{\Q}/(\alg)$. 
\end{lem}

\Pf . This can be easily deduced from (\ref{Ufor}) and 
(\ref{Dfor}) using the following
commutation relations for the operators $\De$, $\del_1$ and $H$:
$$[\del_1,\De]=0,\ [H,\del_1]=-\del_1,\ [H,\De]=-2\De.$$
\ed

The powers of the operator $\De$ are computed in the following lemma.

\begin{lem}\label{Depowerlem}
For $m\ge 0$, $n\ge 0$ and $n_1\ldots,n_k$ such that $n_i>1$ one has
$$\De^{[m]}(x_1^{[n]}x_{n_1}\ldots x_{n_k})=
\sum_{[1,k]=I_1\sqcup I_2\sqcup\ldots\sqcup I_{k-m}}b(I_1)\ldots b(I_{k-m})
x_1^{[n]}x_{d(I_1)}\ldots x_{d(I_{k-m})},
$$
where we use the notation of Theorem \ref{mainthm},
the summation is over all unordered
partitions of the set $[1,k]=\{1,\ldots,k\}$
into the disjoint union of $k-m$ nonempty subsets.
\end{lem}

\Pf . For $m=0$ the assertion is clear, so we can use induction in $m$.
Using the definition of $\De$ we can immediately reduce the induction
step to the following identity:
$$\sum_{[1,k]=I\sqcup J}b(I)b(J){d(I)+d(J)\choose d(I)}=2(k-1)b([1,k])$$
of polynomials in $n_1,\ldots,n_k$, where in the LHS we consider
partitions of $[1,k]$ into the {\it ordered} disjoint union of nonempty subsets $I$ and $J$. 
Equivalently, we have to check that
\begin{equation}\label{combid1}
\sum_{[1,k]=I\sqcup J}P_{|I|-1}(\sum_{i\in I} u_i)P_{|J|-1}(\sum_{j\in J} u_j)=
2(k-1)P_{k-2}(\sum_{i=1}^k u_i)
\end{equation}
where $u_1,\ldots,u_k$ are formal commuting variables;
the sum is of the same kind as before; 
$P_n(u):=u(u-1)\ldots (u-n+1)$ for $n\ge 0$
($P_0=1$).
Note that the formal power series $P_t(u)=\sum_{n\ge 0}P_n(u)t^n/n!=(1+t)^u$
satisfies $P_t(u+v)=P_t(u)P_t(v)$. Hence,
\begin{equation}\label{multeq}
P_n(u+v)=\sum_{i=0}^n {n\choose i} P_i(u)P_{n-i}(v).
\end{equation} 
Using this property we can express both parts of the identity (\ref{combid1}) 
as linear combinations of the products of the form $P_{i_1}(u_1)\ldots P_{i_k}(u_k)$,
where $i_1+\ldots+i_k=k-2$. Thus, (\ref{combid1}) is equivalent to a sequence of identities
obtained by equating the coefficients with each such product. Note that the polynomials
$\wt{P}_n(u)=u^n$ also satisfy (\ref{multeq}). Hence, (\ref{combid1})
is equivalent to a similar identity with $P_n(u)$ replaced by $u^n$:
\begin{equation}\label{combid2}
\sum_{[1,k]=I\sqcup J}(\sum_{i\in I} u_i)^{|I|-1}(\sum_{j\in J} u_j)^{|J|-1}=
2(k-1)(\sum_{i=1}^k u_i)^{k-2}.
\end{equation}
This is the so called Hurwitz identity proved in \cite{H} (see also 
Exer. 5.31b of \cite{St} for a more recent treatment).
\ed

Lemmata \ref{Upowerlem} and \ref{Depowerlem}
give an explicit formula for the action of all the operators
$U^{[m]}$ on monomials in $p_i$'s. Using this formula
we can easily prove our first main theorem.

\subsection{Proof of Theorem \ref{mainthm}}

(i) We start with the obvious vanishing
$$p_1^{m_1}p_2^{m_2}\ldots p_k^{m_k}=0$$
in $\CH(J)_{\Q}/(\alg)$,
where $m_1+2m_2+\ldots+km_k=g$ and $m_1<g$. 
Indeed, this is a class of a cycle of dimension $0$ 
that is homologically equivalent to zero,
hence, it is also algebraically equivalent to zero.
Therefore, for every $d\ge 0$ we have
$$U^d(p_1^{m_1}p_2^{m_2}\ldots p_k^{m_k})=0$$
in $\CH(J)_{\Q}/(\alg)$.
It remains to use (\ref{Ufor}).

\noindent
(ii) Note that for $\sum_{i=1}^k n_i>g$ the relation becomes
trivial, so we can assume that $\sum_{i=1}^k n_i\le g$.
The relations described in (i) have form
$$U^{[d]}(p_1^{[g-\sum_{i=1}^k n_i]}p_{n_1}\ldots p_{n_k})=0.$$
Applying Lemma \ref{Upowerlem} we can rewrite this
as 
\begin{equation}\label{mainid2}
\sum_{j=0}^{d}{k-1-j\choose d-j}
p_1^{[g+j-d-\sum_{i=1}^k n_i]}
\De^{[j]}(p_{n_1}\ldots p_{n_k})=0,
\end{equation}
where we use the convention $p_1^{[n]}=0$ for $n<0$.
Now Lemma \ref{Depowerlem} shows that the obtained identity
is equivalent to the assertion of the theorem.
\ed

\subsection{}

Next we compute $\exp(p_1)*x$, where $x$ is a monomial in $p_i$'s.

\begin{prop}\label{expconvprop}
For $n\ge 0$ and $n_1,\ldots,n_k$ such that $n_i>1$ one has
\begin{align*}
&(-1)^g\exp(p_1)*(p_1^{[n]}p_{n_1}\ldots p_{n_k})=\\
&\sum_{0\le j\le l\le g}(-1)^l {-j+n-g-1+k+\sum_{i=1}^k n_i\choose l-j}
p_1^{[n-l+j]}\De^{[j]}(p_{n_1}\ldots p_{n_k})
\end{align*}
\end{prop}

\Pf . By Lemma \ref{Uconvlem} we have
$$(-1)^g\exp(p_1)*x=\sum_{l=0}^g(-1)^l U^{[l]}(x).$$
Hence, if $x$ is a polynomial in $p_i$'s then using
Lemma \ref{Upowerlem} we get
$$(-1)^g\exp(p_1)*x=\sum_{0\ge j\ge l\le g}(-1)^l \del_1^{l-j}\De^{[j]}
{H-j\choose l-j} x
$$
which immediately gives the required formula.
\ed

We are going to use Proposition \ref{expconvprop} to
compute $S(p_{n_1}\ldots p_{n_k})$, where $n_i>1$. Then
we will deduce the general case of Theorem \ref{fourthm}
from this using our formula for the operators $U^{[n]}$. 

\begin{lem}\label{expconvexplem}
For $n_1,\ldots,n_k$ such that $n_i>1$ one has
\begin{align*}
&(-1)^g\exp(p_1)*(\exp(-p_1)\cdot p_{n_1}\ldots p_{n_k})=\\
&\sum_{j\ge 0,m\ge 0}(-1)^{j+m}{k-j+\sum_{i=1}^k n_i\choose g-m}
p_1^{[m]}\De^{[j]}(p_{n_1}\ldots p_{n_k}).
\end{align*}
\end{lem}

\Pf . We have
$$(-1)^g\exp(p_1)*(\exp(-p_1)\cdot p_{n_1}\ldots p_{n_k})=
(-1)^g\sum_{n\ge 0}(-1)^n\exp(p_1)*(p_1^{[n]}p_{n_1}\ldots p_{n_k}).$$
Using Proposition \ref{expconvprop} this expression can be rewritten
as follows:
\begin{align*}
&\sum_{n\ge 0}
\sum_{0\le j\le l\le g}(-1)^{n+l} {-j+n-g-1+k+\sum_{i=1}^k n_i\choose l-j}
p_1^{[n-l+j]}\De^{[j]}(p_{n_1}\ldots p_{n_k})=\\
&\sum_{j\ge 0}\De^{[j]}(p_{n_1}\ldots p_{n_k})
\cdot\sum_{m\ge 0}(-1)^{j+m}p_1^{[m]}
\sum_{m\le n\le m+g-j}{-j+n-g-1+k+\sum_{i=1}^k n_i\choose n-m}.
\end{align*}
It remains to use the elementary identity
$$\sum_{i=0}^N {x+i\choose i}={x+N+1\choose N}$$
to simplify the above sum.
\ed

\begin{prop}\label{nop1prop}
For $n_1,\ldots,n_k$ such that $n_i>1$ one has
$$S(p_{n_1}\ldots p_{n_k})=
\sum_{j\ge 0}p_1^{[j+g-k-\sum_{i=1}^k n_i]}
\De^{[j]}(p_{n_1}\ldots p_{n_k}),$$
where we use the notation $p^{[m]}=0$ for $m<0$.
\end{prop}

\Pf . Applying formula (\ref{fourexpeq}) for $x=p_{n_1}\ldots p_{n_k}$
and using Lemma \ref{expconvexplem} we get
\begin{equation}\label{nonhomeq}
(-1)^gS[-1]^*(p_{n_1}\ldots p_{n_k})=
\sum_{j\ge 0,m\ge 0}(-1)^{j+m}{k-j+\sum_{i=1}^k n_i\choose g-m}
\exp(-p_1)p_1^{[m]}\De^{[j]}(p_{n_1}\ldots p_{n_k}).
\end{equation}
Now we observe that the LHS belongs to $\CH^{g-k}(J)_{\Q}/(\alg)$.
Indeed, by Lemma \ref{chernlem} we have
$p_{n_1}\ldots p_{n_k}\in\CH^{\sum n_i}_{-k+\sum n_i}/(\alg)$, so
this follows from (\ref{bigreq}).
Note also that in the RHS of (\ref{nonhomeq}) we can restrict
the summation to $(j,m)$ satisfying the inequality
$$g-m\le k-j+\sum_{i=1}^k n_i$$
(otherwise the binomial coefficient vanishes).
For such $(j,m)$ we have
$$p_1^{[m]}\De^{[j]}(p_{n_1}\ldots p_{n_k})\in\CH^d(J)_{\Q}/(\alg)$$
where $d=\sum_{i=1}^k n_i-j+m\ge g-k$.
Thus, the only terms in the RHS belonging to
$\CH^{g-k}(J)_{\Q}/(\alg)$ correspond to
$m=j+g-k-\sum_{i=1}^k n_i$. Therefore, (\ref{nonhomeq}) implies
$$(-1)^gS[-1]^*(p_{n_1}\ldots p_{n_k})=(-1)^{g+k+\sum_{i=1}^k n_i}
\sum_{j\ge 0}p_1^{[j+g-k-\sum_{i=1}^k n_i]}
\De^{[j]}(p_{n_1}\ldots p_{n_k}).$$
It remains to use the formula $[-1]^*p_n=(-1)^{n+1}p_n$.
\ed

\subsection{Proof of Theorem \ref{fourthm}}

We have
$$S(p_1^{[n]}p_{n_1}\ldots p_{n_k})=U^{[n]}(S(p_{n_1}\ldots p_{n_k})).$$
Hence, using Proposition \ref{nop1prop} and Lemma \ref{Upowerlem}
we get
\begin{align*}
&S(p_1^{[n]}p_{n_1}\ldots p_{n_k})=\\
&\sum_{j\ge 0}\sum_{l=0}^n
\del_1^{n-l}\De^{[l]}{H-l\choose n-l}
p_1^{[j+N]}\De^{[j]}(p_{n_1}\ldots p_{n_k}),
\end{align*}
where $N=g-k-\sum_{i=1}^k n_i$.
Applying the definition of $H$ and $\del_1$ we can rewrite this as 
\begin{align*}
&\sum_{j\ge 0}\sum_{l=0}^n {-j-1-l\choose n-l}
p_1^{[j+l-n+N]}\De^{[l]}\De^{[j]}(p_{n_1}\ldots p_{n_k})=\\
&\sum_{j\ge 0}\sum_{l=0}^n {-j-1-l\choose n-l}{j+l\choose l}
p_1^{[j+l-n+N]}\De^{[j+l]}(p_{n_1}\ldots p_{n_k})=\\
&\sum_{m\ge 0}\sum_{l=0}^n {-m-1\choose n-l}{m\choose l}
p_1^{[m-n+N]}\De^{[m]}(p_{n_1}\ldots p_{n_k}).
\end{align*}
But one has
$$\sum_{l=0}^n {-m-1\choose n-l}{m\choose l}=(-1)^n$$
(this can be proved by looking at coefficients with $t^n$
in the identity $(1+t)^{-m-1}(1+t)^m=(1+t)^{-1}$).
Hence,
\begin{equation}\label{foureq2}
S(p_1^{[n]}p_{n_1}\ldots p_{n_k})=(-1)^n\sum_{m\ge 0}
p_1^{[m-n+g-k-\sum_{i=1}^k n_i]}\De^{[m]}(p_{n_1}\ldots p_{n_k}).
\end{equation}
In view of Lemma \ref{Depowerlem} the obtained identity is equivalent
to Theorem \ref{fourthm}.
\ed

\subsection{}\label{corsec}

\noindent
{\it Proof of Corollary \ref{CVGcor}}.
For $d=k-1$ and $\sum_{i=1}^k n_i=g$ 
the relation (\ref{mainid}) gives $p_{g-k+1}=0$.
It remains to observe that for every $m$ such that $g/2+1\le m\le g$
we can choose $(n_1,\ldots,n_{g-m+1})$ with
$n_i\ge 2$ and $\sum_i n_i=g$. 
\ed

\noindent
{\it Proof of Corollary \ref{lincombcor}}.
We want to represent an element
$$p_1^{[g-d-\sum_{i=1}^k n_i]}p_{n_1}\ldots p_{n_k},$$
where $n_i>1$ and $k>d$,
as a linear combination of classes of the form 
$$p_1^{[g-d-\sum_{i=1}^{k'} n'_i]}p_{n'_1}\ldots p_{n'_{k'}}$$
with $k'\le d$ (and $n'_i>1$).
For this we can use a system of relations
\begin{align*}
&D^{[d]}(p_1^{[m+g-\sum_i n_i]}\De^m(p_{n_1}\ldots p_{n_k}))=\\
&\sum_{j=0}^{d}{k-1-m-j\choose d-j}\cdot\frac{1}{j!}
p_1^{[g+m+j-d-\sum_i n_i]}\De^{m+j}(p_{n_1}\ldots p_{n_k})=0
\end{align*}
for $m=0,\ldots,k-2$. 
For $j\le k-1$ let us denote
$$a_j=p_1^{[g+k-1-j-d-\sum_i n_i]}\De^{k-1-j}(p_{n_1}\ldots p_{n_k}),$$
so that $a_j=0$ for $j<0$. Note that $a_0,\ldots,a_{d-1}$ are linear
combinations of classes of the required form. On the other hand, 
the above relations can be rewritten as 
\begin{equation}\label{recrel}
\sum_{j=0}^d{N-j\choose d-j}\cdot\frac{1}{j!}a_{N-j}=0
\end{equation}
for $N\le k-1$. It is clear that these relations allow to express $a_{k-1}$
in terms of the classes $a_0,\ldots,a_{d-1}$.
To find the explicit formula we observe that we can define $a_j$ for all
$j\in\Z$ by imposing the recursive relation (\ref{recrel}) also for $N\ge k$.
Then the generating function
$$F(t)=\sum_{j\ge 0}a_jt^j$$
satisfies the differential equation
$$\sum_{j=0}^d{d\choose j}F^{(j)}(t)=0.$$
The initial terms $a_0,\ldots,a_{d-1}$ determine the solution uniquely:
$$F(t)=\sum_{j=0}^{d-1}a_jt^j\cdot(\sum_{m=0}^{d-1-j}\frac{t^m}{m!})\cdot
\exp(-t).$$
Hence, we obtain
$$a_N=\sum_{j=0}^{d-1}a_j\cdot
\sum_{m=0}^{d-1-j}\frac{(-1)^{N-j-m}}{m!(N-j-m)!}.$$
Using the elementary identity 
$$\sum_{m=0}^{d-1-j}(-1)^m{N-j\choose m}=(-1)^{d-1-j}{N-1-j\choose d-1-j}$$
we derive that
$$a_N=(-1)^{N+d-1}\sum_{j=0}^{d-1}\frac{a_j}{(N-j)!}{N-1-j\choose d-1-j}.$$
Hence, for $N=k-1$ we get
\begin{align*}
&p_1^{[g-d-\sum_i n_i]}p_{n_1}\ldots p_{n_k}=a_{k-1}=\\
&(-1)^{k+d}\sum_{j=0}^{d-1}{k-2-j\choose d-1-j}
p_1^{[g+k-1-j-d-\sum_i n_i]}\De^{[k-1-j]}(p_{n_1}\ldots p_{n_k}),
\end{align*}
which is equivalent to (\ref{niceeq}).
\ed

Now let us consider the relations (\ref{mainid}) for small genera.
Recall that to produce a nontrivial relation the parameters $(d,k,n_1,\ldots,n_k)$ 
should satisfy the inequalities $\sum_{i=1}^k n_i\le g$, $d\le k-1$.
In particular, we should have $d\le k-1\le g/2-1$. Thus, for $g=4$ the only interesting
relation is $p_3=0$ (corresponding to $d=1$, $k=2$, $n_1=n_2=2$).
Here is the complete list of relations equivalent to (\ref{mainid}) for $5\le g\le 10$
(where we omit the relations in codimension $g$):

\noindent $g=5$: 
\begin{equation}
\begin{array}{l}
p_4=0, \ p_2^2=-6p_1p_3;\nonumber
\end{array}
\end{equation}

\noindent $g=6$: 
\begin{equation}
\begin{array}{l}
d=2:\ p_4=0,\nonumber\\ 
d=1:\ p_5=0,\ p_2p_3=0,\ p_1p_4=0,\ p_1p_2^2=-3p_1^2p_3;\nonumber
\end{array}
\end{equation}

\noindent $g=7$:
\begin{equation}
\begin{array}{l}
d=2:\ p_5=0,\ p_2p_3=-5p_1p_4,\nonumber\\ 
d=1:\ p_6=0,\ p_2p_4=0,\ p_3^2=0,\ p_2^3=45p_1^2p_4,\ p_1p_5=0,\
p_1p_2p_3=-5p_1^2p_4,\nonumber\\
p_1^2p_2^2=-2p_1^3p_3;\nonumber
\end{array}
\end{equation}

\noindent $g=8$: 
\begin{equation}
\begin{array}{l}
d=3:\ p_5=0,\nonumber\\ 
d=2:\ p_6=0,\ 3p_3^2=-10p_2p_4,\ p_2^3+18p_1p_2p_3+45p_1^2p_4=0,\ p_1p_5=0,
\nonumber\\
d=1:\ p_7=0,\ p_2p_5=0,\ p_3p_4=0,\ p_2^2p_3=0,\ p_1p_6=0,\ p_1p_2p_4=0,\
p_1p_3^2=0,\nonumber\\
p_1p_2^3=15p_1^3p_4,\ p_1^2p_5=0,\ 3p_1^2p_2p_3=-10p_1^3p_4,\
2p_1^3p_2^2=-3p_1^4p_3;\nonumber
\end{array}
\end{equation}

\noindent $g=9$: 
\begin{equation}
\begin{array}{l}
d=3:\ p_6=0,\ 3p_3^2+10p_2p_4+70p_1p_5=0,\nonumber\\
d=2:\ p_7=0,\ p_2p_5=0,\ p_3p_4=0,\ 2p_2^2p_3+10p_1p_2p_4+3p_1p_3^2=0,\
p_1p_6=0,\nonumber\\
3p_1p_3^2+10p_1p_2p_4+70p_1^2p_5=0,\ p_1p_2^3+9p_1^2p_2p_3+15p_1^3p_4=0,
\nonumber\\
d=1:\ p_8=0,\ p_2p_6=0,\ p_3p_5=0,\ p_4^2=0,\ p_2^2p_4=0,\ p_2p_3^2=0,\
p_2^4=-420p_1^3p_5,\nonumber\\
p_1p_7=0,\ p_1p_2p_5=0,\ p_1p_3p_4=0,\ p_1p_2^2p_3=35p_1^3p_5,\
p_1^2p_2p_4=-5p_1^3p_5,\nonumber\\
3p_1^2p_3^2=-20p_1^3p_5,\ 2p_1^2p_2^3=15p_1^4p_4,\ 2p_1^3p_2p_3=-5p_1^4p_4,\ 
5p_1^4p_2^2=-6p_1^5p_3;  \nonumber
\end{array}
\end{equation}

\noindent $g=10$:
\begin{equation}
\begin{array}{l}
d=4:\ p_6=0,\nonumber\\
d=3:\ p_7=0,\ p_2p_5=0,\ p_3p_4=0,\
p_2^2p_3+3p_1p_3^2+10p_1p_2p_4+35p_1^2p_5=0,\ p_1p_6=0,\nonumber\\
d=2:\ p_8=0,\ p_2p_6=0,\ p_3p_5=0,\ p_4^2=0,\ p_2^2p_4=0,\ p_2p_3^2=0,
\nonumber\\
p_2^4+24p_1p_2^2p_3+18p_1^2p_3^2+60p_1^2p_2p_4=0,\ p_1p_7=0,\ p_1p_3p_4=0,\
p_1p_2p_5=0,\nonumber\\
p_1p_2^2p_3+3p_1^2p_3^2+10p_1^2p_2p_4+35p_1^3p_5=0,\ 
p_1^2p_6=0,\ 2p_1^2p_2^3+12p_1^3p_2p_3+15p_1^4p_4=0,\nonumber\\
d=1:\ p_9=0,\ p_2p_7=0,\ p_3p_6=0,\ p_4p_5=0,\ p_2^2p_5=0,\ p_2p_3p_4=0,\
p_3^3=0,\ p_2^3p_3=0,\nonumber\\
p_1p_8=0,\ p_1p_2p_6=0,\ p_1p_3p_5=0,\ p_1p_4^2=0,\
p_1p_2^2p_4=0,\ p_1p_2p_3^2=0,\nonumber\\
p_1p_2^4=-105p_1^4p_5,\ p_1^2p_7=0,\ p_1^2p_2p_5=0,\ p_1^2p_3p_4=0,\ 2p_1^2p_2^2p_3=35p_1^4p_5,
\nonumber\\
p_1^3p_6=0,\ 2p_1^3p_2^3=9p_1^5p_4,\ 4p_1^3p_2p_4=-15p_1^4p_5,\
p_1^3p_3^2=-5p_1^4p_5,\nonumber\\
p_1^4p_2p_3=-2p_1^5p_4,\ p_1^5p_2^2=-p_1^6p_3.\nonumber
\end{array}
\end{equation}

\section{Formal consequences of relations}

Throughout this section $x_1,x_2,\ldots$ are formal
commuting variables. We equip the ring of polynomials
$\Q[x_1,x_2,\ldots]$ with grading by setting $\deg(x_i)=i$.

\subsection{} The following lemma will allow us to express the relations
in $R^{\Jac}_g$ in the form asserted in Theorem \ref{fourthm}(i).

\begin{lem}\label{Dsurlem} 
Let $V_n\sub V=\Q[x_1,x_2,\ldots]$ denote the component
of degree $n$. The linear map $D:V_{n+1}\to V_n$ induced by the differential
operator $D$ is surjective for $n>g$.
The image of the map $D:V_{g+1}\ra V_g$ is the codimension $1$ subspace
spanned by all monomials of degree $g$ except for $x_1^g$.
\end{lem}

\Pf . Let us consider the operator $Dx_1:V_n\to V_n$.
Then for every monomial $F=x_1^{n-\sum_i n_i}x_{n_1}\ldots x_{n_k}$ in $V_n$,
where $n_i>1$, one has
$$Dx_1F\equiv (n+k-g)F\ \mod (x_1^{n-\sum_i n_i+1})\cap V_n,$$
where $(x_1^m)$ denotes the ideal generated by $x_1^m$. This immediately
implies that $Dx_1$ surjects onto $V_n$ for $n>g$ and that the image of
$Dx_1|_{V_g}$ is the subspace spanned by all monomials of degree $g$ except
for $x_1^g$. It remains to observe that for an arbitrary monomial
$F$ of degree $g+1$, $D(F)$ lies in this subspace. 
\ed

\begin{lem}\label{sl2actlem}
The operators $e(F)=x_1\cdot F$,
$f=-D$ and $h=-g+\sum_{n\ge 1}(n+1)x_n\del_n$
define a representation of $\ssl_2$ on $\Q[x_1,x_2,\ldots]$. 
\end{lem}

The proof is straightforward and is left for the reader.

\begin{lem}\label{sl2invlem} 
Let $V$ be an $\ssl_2$-module. 
Then the subspace $\cap_{n\ge 0}f^n(V)\sub V$ is $\ssl_2$-invariant.
\end{lem}

\Pf . We have $hf^n(V)\sub f^n(V)$. Now an easy induction in $n$ shows that
$ef^n(V)\sub f^{n-1}(V)$. Hence, $\cap_{n\ge 0}f^n(V)$ is invariant with respect to $e$.
\ed

Let $S$ be the operator on $\Q[x_1,x_2,\ldots]$ defined by
\begin{equation}\label{foureq3}
S(x_1^{[n]}x_{n_1}\ldots x_{n_k})=(-1)^n\sum_{m\ge 0}
x_1^{[m-n+g-k-\sum_{i=1}^k n_i]}\De^{[m]}(x_{n_1}\ldots x_{n_k}),
\end{equation}
where $n_i>1$. Of course, this is just the formula (\ref{foureq2}) with
$p_i$ replaced by $x_i$ (recall that we use the convention $x_1^{[d]}=0$
for $d<0$).

\begin{lem}\label{foursqlem} 
For arbitrary $n\ge 0$, $n_1,\ldots, n_k$ such that $n_i>1$ and
$n+k+\sum_{i=1}^k n_i\le g$ one has
$$S^2(x_1^{[n]}x_{n_1}\ldots x_{n_k})=(-1)^{g+k+\sum_i n_i}
x_1^{[n]}x_{n_1}\ldots x_{n_k}.$$
\end{lem}

\Pf .
The condition $n+k+\sum_i n_i\le g$ implies
that in the sum defining $S(x_1^{[n]}x_{n_1}\ldots x_{n_k})$
all terms contain nonnegative powers of $x_1$, i.e., 
no terms are omitted due to our
convention that $x_1^{[d]}=0$ for $d<0$. Therefore, we can apply the formula for $S$
again and find that
\begin{align*}
&(-1)^{g+k+\sum_i n_i}S^2(x_1^{[n]}x_{n_1}\ldots x_{n_k})=
\sum_{m\ge 0,l\ge 0}(-1)^m {m+l\choose m}x_1^{[l+m+n]}\De^{[l+m]}(x_{n_1}\ldots x_{n_k})=\\
&\sum_{d\ge 0}\sum_{m+l=d}(-1)^m {d\choose m}x_1^{[d+n]}\De^{[d]}(x_{n_1}\ldots x_{n_k}).
\end{align*}
It remains to use the fact that $\sum_{m=0}^d (-1)^m {d\choose m}=0$ for $d\ge 1$.
\ed

\subsection{Proof of Theorem \ref{formalthm}(i),(iii),(iv)}

Part (i) is an immediate consequence of Lemma \ref{Dsurlem}.

By Lemma \ref{sl2actlem} we have an action of $\ssl_2$ on 
$V=\Q[x_1,x_2,\ldots]$. Applying Lemma \ref{sl2invlem} and part (i) 
we derive that the subspace $I=I_g=\cap_{n\ge 0}f^n(V)\sub V$ is $\ssl_2$-invariant.

Now let us consider the operator $S:V\ra V$ defined by (\ref{foureq3}).
It is easy to check that 
\begin{equation}\label{Seeq}
Se=-fS.
\end{equation}
Indeed, the formula (\ref{foureq3}) was constructed in such a way that
$$S(x_1^{[n]}x_{n_1}\ldots x_{n_k})=D^{[n]}S(x_{n_1}\ldots x_{n_k}),$$ 
where $n_i>1$,
Since $f=-D$ this implies the above relation.
It is also straightforward to check that
\begin{equation}\label{Sheq}
Sh=-hS.
\end{equation}
Next we claim that
\begin{equation}\label{Sfeq}
Sf(F)+eS(F)\in I
\end{equation}
for every $F\in V$.
Using the relation (\ref{Seeq}) one can easily check that
$$(Sf+eS)e=-f(Sf+eS).$$ 
Hence, the claim would follow if we could prove that
$$(Sf+eS)(x_{n_1}\ldots x_{n_k})\in I,$$
where $n_i>1$. We have
\begin{align*}
&(Sf+eS)(x_{n_1}\ldots x_{n_k})=x_1S(x_{n_1}\ldots x_{n_k})-
S\De(x_{n_1}\ldots x_{n_k})=\\
&\sum_{m\ge 0}(g+m+1-k-\sum n_i)x_1^{[g+m+1-k-\sum n_i]}
\De^{[m]}(x_{n_1}\ldots x_{n_k})-\\
&\sum_{m\ge 0}m x_1^{[g+m+1-k-\sum n_i]}\De^{[m]}(x_{n_1}\ldots x_{n_k})=\\
&(g+1-k-\sum n_i)
\sum_{m\ge 0}x_1^{[g+m+1-k-\sum n_i]}\De^{[m]}(x_{n_1}\ldots x_{n_k}).
\end{align*}
If $\sum n_i>g$ then the latter sum is zero since nonzero terms
correspond to $m\le k-1$ and $g+m+1-k-\sum n_i\ge 0$.
Otherwise, it is proportional to $D^{[k-1]}(x_1^{g-\sum n_i}x_{n_1}\ldots x_{n_k})$
(see the proof of Theorem \ref{mainthm}(ii)).
Hence, it belongs to $I$ which proves our claim.

Next, we observe that if $n+\sum_{i=1}^k n_i>g$ (where $n_i>1$) then
$$S(x_1^{[n]}x_{n_1}\ldots x_{n_k})=(-1)^n\sum_{m>k}
x_1^{[m-n+g-k-\sum_{i=1}^k n_i]}\De^{[m]}(x_{n_1}\ldots x_{n_k})=0$$
since $\De^{[m]}(x_{n_1}\ldots x_{n_k})=0$ for $m>k$.
Hence $S(V_m)=0$ for $m>g$, where $V_m$ is the component of degree $m$ in $V$. 
Using (\ref{Sfeq}) and Lemma \ref{Dsurlem}
we conclude that $S(I)\sub I$.

Finally, we want to prove that
\begin{equation}\label{S2eq}
S^2(F)=(-1)^g[-1]^*(F)\mod I
\end{equation}
for every $F\in V$, where
$[-1]^*x_1^{[n]}x_{n_1}\ldots x_{n_k}=(-1)^{k+\sum n_i}
x_1^{[n]}x_{n_1}\ldots x_{n_k}$.
The relations (\ref{Seeq}), (\ref{Sheq}) and (\ref{Sfeq}) imply that
the operator $S^2-(-1)^g[-1]^*$ induces an endomorphism of the finite-dimensional
$\ssl_2$-module $V/I$. Therefore, it is enough to prove that (\ref{S2eq}) holds
for every $F\in V$ such that $hF=mF$ with $m\le 0$. Thus, it suffices to prove
that
$$S^2(x_1^{[n]}x_{n_1}\ldots x_{n_k})=(-1)^{g+k+\sum n_i}
x_1^{[n]}x_{n_1}\ldots x_{n_k},$$
where $n_i>1$, provided that $-g+2n+k+\sum n_i\le 0$.
But in this case we also have $-g+n+k+\sum n_i\le 0$, 
so we are done by Lemma \ref{foursqlem}.
\ed

\subsection{}
For $k\ge 0$, $m\ge 0$, $k+m\ge 2$ consider the differential operators
$$\De_{k,m}=\frac{1}{k!}
\sum_{n_1,\ldots,n_k\ge 2}\frac{(n_1+\ldots+n_k+m)!}{n_1!\ldots n_k!m!}
x_{n_1+\ldots+n_k+m-1}\del_{n_1}\ldots\del_{n_k}.$$
Note that $\De_{0,m}=x_{m-1}$ for $m\ge 2$ and that 
$\De_{2,0}=\De$.
These operators appear when we want to write the operators $D_k$ 
considered in Theorem \ref{diffthm} as polynomials in $\del_1$:
for every $k\ge 2$ one has
\begin{equation}\label{DDeeq}
D_k=\sum_{m=0}^k\De_{k-m,m}\del_1^m.
\end{equation}
Hence, we have
$$D=D_2-g\del_1=\De_{2,0}+\De_{1,1}\del_1+p_1\del_1^2-g\del_1.$$ 

\begin{lem}\label{Dcomrellem} 
(i) The following commutation relations hold:
$$[\De_{k,m},\De_{k',m'}]=(km'-k'm)\frac{(k+k'-1)!(m+m'-1)!}{k!k'!m!m'!}\De_{k+k'-1,m+m'-1}$$
for $m+m'\ge 1$ and
$$[\De_{k,0},\De_{k',0}]=0.$$
In particular,
$$[\De,\De_{k,m}]=(k+1)\De_{k+1,m-1}$$
for $m\ge 1$ and 
$$[\De,\De_{k,0}]=0.$$

\noindent
(ii) One has $[D_k,D_{k'}]=0$ for $k,k'\ge 2$ and
$[D_k,\del_1]=0$ for $k\ge 3$.
\end{lem}

We leave the proof for the reader.

Let us denote by $V'\sub V$ the subspace spanned by all the monomials
$x_1^nx_{n_1}\ldots x_{n_k}$, where $n_i>1$, such that $n+k+\sum_{i=1}^k n_i\le g$.
Let also $J\sub I$ denote the span of the subspaces $V_m$ for $m>g$ together with all the elements
of the form
$$
D^{[d]}(x_1^{[g-\sum_i n_i]}x_{n_1}\ldots x_{n_k})=
\sum_{j=0}^{d}{k-1-j\choose d-j}
x_1^{[j+g-d-\sum_i n_i]}
\De^{[j]}(x_{n_1}\ldots x_{n_k}),
$$
where $k\ge 1$, $n_i>1$, $d\ge 0$ and $d+\sum_i n_i\le g$.

\begin{lem}\label{importantlem} 
(i) $V=J+V'$;

\noindent (ii) $I\cap V'=S(I)=S(I\cap V')$;

\noindent (iii) $I=J+S(I\cap V')$.
\end{lem} 

\Pf .(i) This follows essentially from the proof of Corollary \ref{lincombcor}.
Indeed, it suffices to show that every monomial
$x_1^nx_{n_1}\ldots x_{n_k}$ such that $n_i>1$, $n+k+\sum n_i>g$,
can be expressed modulo $J$ in terms of monomials 
$x_1^{n'}x_{n'_1}\ldots x_{n'_{k'}}$ with $n'+k'+\sum n'_i<n+k+\sum n_i$.
If $n+\sum n_i>g$ then our monomial belongs to $J$, so we can assume that
$d:=g-n-\sum n_i\ge 0$. Our assumption implies that $d<k$, hence
the element $D^{[d]}(x_1^{[g-\sum_i n_i]}x_{n_1}\ldots x_{n_k})\in J$
has a nonzero coefficient with $x_1^{[n]}x_{n_1}\ldots x_{n_k}$. Therefore, it
gives the required expression.

\noindent (ii) The inclusion $S(V)\sub V'$ is immediate from the definition of $S$.
Since $S$ preserves $I$ we obtain $S(I\cap V')\sub S(I)\sub I\cap V'$. On the other hand, Lemma
\ref{foursqlem} implies that $S^2|_{V'}=(-1)^g[-1]^*$. Hence, $I\cap V'\sub S(I\cap V')$.

\noindent (iii) Combining (i) and (ii) we get 
$I=J+I\cap V'=J+S(I\cap V')$.
\ed

\begin{rem} In fact, it is easy to show that $J=\ker(S)$. Indeed, the same argument
as in the proof below shows that $S(F)=0$ for $F\in J$. On the other hand,
$\ker(S)\cap V'=0$ by Lemma \ref{foursqlem}. Hence, $J=\ker(S)$. 
It follows that we have direct sum decompositions
$V=J\oplus V'$, $I=J\oplus S(I\cap V')$.
\end{rem}

\subsection{Proof of Theorem \ref{formalthm}(ii) and of Theorem \ref{diffthm}}

First, we observe that by Lemma \ref{Dcomrellem}(ii) 
the operators $D_k$ for $k\ge 3$
commute with $D=D_2-g\del_1$, hence they preserve $I_g$ and induce operators on
$R^{\Jac}_g$.

Next, we claim that the identity
\begin{equation}\label{Sxmeq}
S(x_m(S(F)))=(-1)^gD_{m+1}[-1]^*F
\end{equation}
holds for every $F\in V'$, $m\ge 2$.
It suffices to take $F=x_1^{[n]}x_{n_1}\ldots x_{n_k}$, where $n_i>1$, $n+k+\sum_i n_i\le g$.
Applying the definition of $S$ we can write 
\begin{align*}
&Sx_mS(x_1^{[n]}x_{n_1}\ldots x_{n_k})=
(-1)^n\sum_{i\ge 0} Sx_1^{[i-n+g-k-\sum_s n_s]}\De^{[i]}(x_{n_1}\ldots x_{n_k})=\\
&(-1)^{g-k-\sum_s n_s}
\sum_{i,j\ge 0}(-1)^ix_1^{[i+j+n-m-1]}\De^{[j]}x_m\De^{[i]}(x_{n_1}\ldots x_{n_k}),
\end{align*}
where we used the assumption $n+k+\sum_i n_i\le g$ to be able to apply the definition
of $S$ the second time.
Now we observe that for every $N\ge 0$ we have the equality of operators
$$\sum_{i+j=N}(-1)^i\De^{[j]}x_m\De^{[i]}=\frac{1}{N!}(\ad\De)^N(x_m).$$
Recall that $x_m=\De_{0,m+1}$, so using Lemma \ref{Dcomrellem}(i) we find
$$(\ad\De)^N(x_m)=N!\De_{N,m+1-N}$$
for $N\le m+1$ and $(\ad\De)^N(x_m)=0$ for $N>m+1$.
Therefore, we obtain
$$Sx_mS(x_1^{[n]}x_{n_1}\ldots x_{n_k})=
(-1)^{g-k-\sum_i n_i}\sum_{N=0}^{m+1}x_1^{[N+n-m-1]}\De_{N,m+1-N}(x_{n_1}\ldots x_{n_k}).$$
Comparing the result with (\ref{DDeeq}) we observe that it is equal to
$$(-1)^{g-k-\sum_i n_i}D_{m+1}(x_1^{[n]}x_{n_1}\ldots x_{n_k})$$
which proves (\ref{Sxmeq}).

It follows that $Sx_mS(I\cap V')=D_{m+1}(I\cap V')\sub I$, hence
\begin{equation}\label{incl1}
x_mS(I\cap V')\sub I.
\end{equation}
Now we are going to show that $S(x_mF)=0$ for $F\in J$. If $F$ has degree $>g$ then
this is clear, so it is enough to check that
$$S\left(x_m\cdot\sum_{j=0}^{d}{k-1-j\choose d-j}
x_1^{[j+g-d-\sum_i n_i]}
\De^{[j]}(x_{n_1}\ldots x_{n_k})\right)=0,$$
where $k\ge 1$, $n_i>1$, $d\ge 0$ and $d+\sum_i n_i\le g$.

Since $S(x_mJ)=0$ we derive that 
\begin{equation}\label{incl2}
x_mJ\sub I. 
\end{equation}
Recall that $I=J+S(I\cap V')$ by Lemma \ref{importantlem}(iii), so
the inclusions (\ref{incl1}) and (\ref{incl2}) imply
that $x_mI\sub I$ for every $m\ge 2$. Since we also know that $I$ is invariant under $e$,
it follows that $I$ is an ideal in $\Q[x_1,x_2,\ldots]$.

Now the decomposition $V=V'+I$ implies 
that equation (\ref{Sxmeq}) holds for all $F\in R^{\Jac}_g$.
This proves (\ref{diffpk}).
\ed

\end{document}